 \documentclass[final,1p,times]{elsarticle}

\usepackage{lineno,hyperref}
\modulolinenumbers[5]

\journal{Journal of \LaTeX\ Templates}

\usepackage{amsmath}
\usepackage{breqn}

\usepackage{numcompress}

\usepackage{epsfig}

\usepackage{amssymb}
\usepackage{float}

\journal{Journal of Mathematical Analysis and Applications}

\begin{document}
%\linenumbers.

\begin{frontmatter}



\title{Fractional Powers of the Differentiation and Integration Operators and its Application in Accelerator Physics and Technology}


\author{M. Behtouei$^{a}$, L. Faillace$^{a}$, L. Palumbo$^{b, c}$, B. Spataro$^{a}$, A. Variola$^{a}$ and M. Migliorati$^{b, c}$\\\vspace{6pt}}

\address{$^{a}${INFN, Laboratori Nazionali di Frascati, P.O. Box 13, I-00044 Frascati, Italy};\\
 $^{b}${Dipartimento di Scienze di Base e Applicate per l'Ingegneria (SBAI), Sapienza University of Rome, Rome, Italy};\\
  $^{c}${INFN/Roma1, Istituto Nazionale di Fisica Nucleare, Piazzale Aldo Moro, 2, 00185, Rome, Italy }}

\begin{abstract}
In this paper we present a solution to a fractional integral of the order 3/2 with the use of a novel method. The integral arises during solving the Biot-Savart equation to find the exact analytical solution for the magnetic field components of a solenoid. We solved the integral by cutting the branch line in order to have an analytic function inside the integral instead of multi-valued operation. 
\end{abstract}




\end{frontmatter}


\section{Introduction}

Fractional calculus studies the fractional powers of the differentiation and integration operators both in the real and complex number domains. The concept was initially brought to the literature by Gottfried Wilhelm Leibniz \cite{ref1} and then developed by Niels Henrik Abel \cite{ref2}, Riemann-Liouville, Riesz, Caputo (see \cite{ref3}), Hadamard \cite{ref4}, Atangana-Baleanu \cite{ref5} and Kober \cite{ref6}. The fractional derivatives and integrals appear in mathematics, in different branch of physics and in engineering. Physicists usually avoid to confront with fractional integrals and derivatives and usually try to find other ways to solve their problems like using the elliptical integral. The idea recently was introduced to quantum mechanics by Nick Laskin \cite{ref6} under the name "fractional quantum mechanics" in which the Brownian-like quantum paths substitute the Lévy-like ones in the Feynman path integral. In this paper we report the main equations of fractional calculus, investigate the limitation of each ones, comparing with other methods, and finally propose a novel way to solve some important of fractional integrals in accelerator physics and engineering. 
It should be noted that the integer derivative of a function f(x)  is a local property at a point x. On the other hand only for non-integer power derivatives, the fractional derivative of a function f(x) at the point x depends only on values of f close to x. This means that the boundary conditions in the theory should be considered, in general, by involving information on the function further out.

In the next section we review the current theorem of fractional derivatives and integrals.

\section{Preliminaries: Review on Fractional Calculus }

In this section we report the several known forms of the fractional integral theorems in mathematics for solving physics and engineering problems.

\subsection{Riemann-Liouville Fractional Integral}

Definition: Let Re $\alpha>0$ and f be piecewise continuous and integrable on (0, $\infty$). Then for $z>0$ we define 

\begin{equation}\label{1}
D_{z-z_0}^{-\alpha} f(z)=\frac{1}{\Gamma(\alpha)} \int_{z_0}^z \frac{f(z)}{(z-z_0)^{-\alpha+1} }dz
\end{equation}
when $z=z_0$, this is called Riemann-Liouville fractional integral of the function $f$ of order $\alpha$.

\subsection{Weyl fractional derivatives}

The Weyl fractional derivatives is used when $z_0$ takes on a singular value from $-\infty$ to $\infty$, and it can be expressed as,

\begin{equation}\label{2}
D_{z^{-\infty}}^{-\alpha} f(z)=\frac{(-1)^{-\alpha}}{\Gamma(-\alpha)} \int_{z_0}^\infty \frac{f(z)^n}{(z-z_0)^{-\alpha-1} }dz
\end{equation}
and 

\begin{equation}\label{3}
D_{z^{+\infty}}^{-\alpha} f(z)=\frac{1}{\Gamma(-\alpha)} \int_{-\infty}^{z_0} \frac{f(z)^n}{(z-z_0)^{-\alpha-1} }dz.
\end{equation}

where Re $\alpha>0$.
\subsection{Caputo fractional derivative}

The Caputo fractional derivative \cite{ref7} is used in order to solve the differential equations without defining the fractional order initial conditions. Caputo's definition is as follows.

\begin{equation}
\small{D_{z-z_0}^{-\alpha} f(z)=\frac{1}{\Gamma(n-\alpha)} \int_{z_0}^z \frac{f(z)^n}{(z-z_0)^{\alpha+1+n} }dz, \ \ n-1<\alpha<n.}
\end{equation}

\subsection{Hadamard Fractional Integral}

The Hadamard fractional integral is introduced by Jacques Hadamard and is given by the following formula,

\begin{equation}
_a D_t^{-\alpha} f(z)=\frac{1}{\Gamma(\alpha)} \int_{a}^t (log\frac{t}{\tau})^{\alpha-1} f(\tau)\frac{d\tau}{\tau}\ \ \ \ t>a
\end{equation}

This is based on the generalization of the integral

\begin{equation}
\int_a^t \frac{d\tau_1}{\tau_1} \int_a^{\tau_1}\frac{d\tau_2}{\tau_2} . . .  \int_a^{\tau_{\alpha-1}}\frac{d\tau_\alpha}{\tau_\alpha}d\tau_\alpha=\frac{1}{\Gamma(\alpha)} \int_{a}^t (log\frac{t}{\tau})^{\alpha-1} f(\tau)\frac{d\tau}{\tau}\ \ \ \ t>a,\ \ \alpha>0
\end{equation}

to obtain the above equation, the n-fold integral of the form below is used 

\begin{equation}
_a^\rho D_x^\alpha f(x)=\int_a^x \tau_1^\rho d\tau_1 \int_a^{\tau_1} \tau_2^\rho d\tau_2 ... \int_a^{\tau_{n-1}} \tau_{n}^\rho  f(\tau_n) d\tau_n
\end{equation}

\subsection{Generalized Fractional Integration Operator }

The author  \cite{ref8} obtained a generalized fractional integration operator which  bounded in the Lebesgue measurable space. The procedure is as follow:

The Lebesgue measurable functions f on [a, b] in the space  $X_c^p$ (a,b) (c $\in {\rm I\!R}$, $1\leq p \leq \infty$) for which $||f_{X_c^p}|| < \infty$, where the norm is defined by

\begin{equation}
||f||_{X_c^p}=\frac{1}{(\int_a^b |t^c f(t)|^p \frac{dt}{t})^p} < \infty \ \ \  (c \in {\rm I\!R}, 1\leq p < \infty)
\end{equation}

and for the case of p when is equal $\infty$ we have

\begin{equation}
||f||_{X_c^\infty}=\mathrm {ess\ sup}_{a \leq t \leq b} [|t^c f(t)|]\ \ \ (c \in {\rm I\!R})
\end{equation}

where \textquotedblleft{ess sup}\textquotedblright  denotes for the \textquotedblleft{essential supremum}\textquotedblright  of the function f, representing that value where f is larger or equal than the function values everywhere, when ignoring what the function does at a set of points of \textquotedblleft{measure zero}\textquotedblright  (measure zero is a set of points capable of being enclosed in intervals whose total length is arbitrarily small). 

By using Dirichlet technique for n-fold integral, the fractional integral of $_a^\rho D_x^\alpha$ yields \cite{ref8},

\begin{equation}
_a^\rho D_x^\alpha f(x)=\frac{(\rho+1)^{1-\alpha}}{\Gamma(\alpha)} \int_a^x (x^{\rho+1}-\tau^{\rho+1})^{\alpha-1} \tau^\rho f(\tau) d\tau
\end{equation}

where $\alpha$ and $\rho\neq -1$ are real numbers.

As an additional information, the Dirichlet technique \cite{ref9} is given:

\begin{equation}
\int_a^x \tau_1^\rho d\tau_1 \int_a^{\tau_1} \tau^\rho  f(\tau) d\tau=\int_a^x \tau^\rho f(\tau) d\tau \int_\tau^{x} \tau_1^\rho   d\tau 
\end{equation}

\begin{equation}
=\frac{1}{\rho+1} \int_a^x (x^{\rho+1}-\tau^{\rho+1})\tau^\rho f(\tau) d\tau
\end{equation}

\subsection{Fractional Atangana-Baleanu derivative}

Atangana and Baleanu proposed a new fractional derivative \cite{ref10} with non-local and no-singular kernel using the Mittag-Leffler function. They started with the fractional
ordinary differential equation 

\begin{equation}
\frac{d^\alpha}{dx^\alpha}=ay\ \ \ 0<\alpha<1
\end{equation}
after some manipulation and using Caputo-Fabrizio derivative they obtained an expression and solved the problem of non-locality. 

\begin{equation}
D_t^\alpha=\frac{M(\alpha)}{1-\alpha} \Sigma_{k=0}^\infty \frac{(-a)^k}{\Gamma (\alpha k+1)} \int_b^t \frac{df(y)}{dy} ((t-y))^{\alpha k} dy.
\end{equation}

where $a=\alpha(1-\alpha)^{-1}$

\subsection{Fractional Riesz derivative}

The Riesz derivative of function u(x, t) with respect to x is defined by \cite{ref11,ref12}

\begin{equation}
\frac{\partial^\alpha u(x,t)}{\partial |x|^\alpha}=-\frac{1}{2} \sec(\frac{\pi \alpha}{2}) [_{RL} D_{-\infty,x}^\alpha+ _{RL} D_{X,+\infty}^\alpha] u(x,t)
\end{equation}

where $D_{-\infty,x}^\alpha$ and  $_{RL} D_{X,+\infty}^\alpha$ are the left and right Riemann-Liouville derivatives and $\Gamma$ is the Euler's Gamma function.

\subsection{The Cauchy integral theorem}

Let $D \subseteq  \mathbb{C}$ be an open set. We shall say that D has a piecewise $C^1$-boundary if the boundary of D (in $\mathbb{C}$) is a closed piecewise $C^1$-contour $\Gamma$ such that each point of $\Gamma$ is also a boundary point of $\mathbb{C}$ $\setminus$ $\bar D$.

 Let $D \subseteq  \mathbb{C}$  be a bounded open set with piecewise $C^1$-boundary, let E be a Banach space, and let f : $\bar D \rightarrow E$ be a continuous function which is holomorphic in D. Then \cite{ref13}

\begin{equation}
f(z_0)=\frac{1}{2\pi i}\int_{\partial D} \frac{f(z)}{z-z_0} dz,\ \ \ z_0 \in D.
\end{equation}

Any holomorphic function with values in a Banach space is infinitely times complexly differentiable. In particular, it is of class $C^\infty$. Moreover, if D and f are as in the above theorem and if we denote by $f^{(n)}$ the n-th complex derivative of f in D, then \cite{ref13}

\begin{equation}
f^{(n)}(z_0)=\frac{n!}{2\pi i}\int_{\partial D} \frac{f(z)}{(z-z_0)^{n+1}} dz,\ \ \ z_0 \in D.
\end{equation}

\section{Statement of the Problem}

Starting from the Biot-Savart law, the axial and radial magnetic field components for a coil of negligible thickness with a stationary electric current are given by \cite{ref14}:

\begin{equation}
B_r=\frac{\mu_0  I  M_z}{4 \sqrt{2}\pi R  } \bigg(\frac{\xi}{\eta}\bigg)^{3/2}\ I_2(\xi)
\end{equation}

\begin{equation}
B_z=\frac{\mu_0  I  }{4\sqrt{2} \pi R } \bigg(\frac{\xi}{\eta}\bigg)^{3/2}\ ( I_1(\xi)-  \eta\ I_2(\xi))
\end{equation}

where $\eta=\frac{r}{R}$, $M_z=\frac{z}{R}$ 

\begin{equation}\label{18}
I_1(\xi)=\int_{0}^{\pi}   \frac{ d\psi}{[1-\xi \ \cos(\psi)]^{3/2}}
\end{equation}

\begin{equation}\label{19}
I_2(\xi)=\int_{0}^{\pi}   \frac{ \cos(\psi)}{[1-\xi \ \cos(\psi)]^{3/2}}\ d\psi
\end{equation}

and

\begin{equation}
\xi (R, z, \eta)=\frac{2\eta }{1+\eta^2+M_z^2}
\end{equation}

and we have used the following notation:

$B_z$ is the magnetic field component in the direction of  the coil axis. 

$B_r$ is the radial magnetic field component.

$I$ is the current in the wire.

$R$ is the radius of the current loop.

z is the distance, on axis, from the center of the current loop to the field measurement point.

$r$ is the radial distance from the axis of the current loop to the field measurement point

$\theta$ denotes for the angle of the current element

$\gamma$ stands for the angle of the observer where the magnetic field components are to be calculated

$\psi=\gamma-\theta$

As it can be observed the fractional integrals (\ref{18}) and (\ref{19}) can not be solved by Riemann-Liouville Fractional Integral. The reason is that, the denominator in this method becomes a power of 3/2 when one choose $\alpha=-1/2$ and this is not allowed by the theorem which implies $\alpha>0$ (see Eq. (\ref{1})). On the other hand, our fractional integrals can not be solved by Weyl fractional theorem ( Eqs. (\ref{2}, \ref{3})), because the $\alpha$ should be equal -5/2 to turn the denominator's power to 3/2 and this is not allowed by the limitation of this theorem ($\alpha>0$). Let us now see if it is possible that the integrals (\ref{18}) and (\ref{19})  can be solved by Caputo fractional theorem. For $n=1$, $\alpha=-1/2$ we can cover the integral's denominators (to arrive at power of 3/2) but as $n-1<\alpha<n$ and this means $\alpha>0$ and this is not allowed either. We face the same limitation considering other theorems mentioned in the previous section. 

In order to solve the integral we use Cauchy's Integral Formula with some modifications. Recalling Cauchy's Integral Formula

\begin{equation}
D^n f(z)=\frac{\Gamma(\alpha +1)}{2\pi i} \oint_C \frac{f(z)}{(z-z_0)^{n+1} }dz
\end{equation} 

we observe that n is an integer number, and this implies that there is one or more singularities in the function. On the other hand, when we have n as a non-integer, the singularity turn to the branch lines and f become a non-local property. It should be noted that integer derivative of a function f is a local property at a point z. We have already observed that only for non-integer power derivatives, the fractional derivative of a function f at the point z depends only on values of f very near z. In order to use this theorem for non-integer n, one should change the multi-valued operation (function) in Eqs. (\ref{18}) and (\ref{19})  and turn it into the analytic function. This procedure is called Branch Cut. By branch cut, our multi-valued function becomes an analytic function with a local property and the branch point turns to be a singularity point. Now we can use the Cauchy's Integral Formula for solving the fractional integral. We call this theorem Fractional Cauchy-like Integral Formula:

\begin{equation}
^{\gamma(z,z^+)} D_{z-z_0}^\alpha f(z)=\frac{ sin (\pi \alpha) \ \Gamma(\alpha +1)}{\pi } \int_{z_0}^z \frac{f(z)}{(z-z_0)^{\alpha+1}}dz.
\end{equation}

In the next section we will first see how we obtain the above equation and then we will apply that to solve our fractional integrals.

\section{Fractional Cauchy-like Integral Formula and the Final Solution}

It should be noted that this method has been studied by the authors of \cite{ref15} and \cite{ref16}. The procedure is as follows:

Let us recall the Cauchy's integral formula

\begin{equation}\label{23}
D^n f(z)=\frac{\Gamma(n +1)}{2\pi i} \oint_C \frac{f(z)}{(z-z_0)^{n+1} }dz.
\end{equation}

Let the contour of integration be $\gamma(z_0,z^+)$. The branch line for $(z-z_0)^{-\alpha - 1}$ starts from the position z and ends at the fixed point $z_0$. The above equation is equivalent to the Riemann-Liouville fractional integral when Re $(\alpha) <0$. We divide the contour $\gamma(z_0,z^+)$  into three contours (see Fig. 1),

\begin{equation}
\gamma(z_0,z^+)=\gamma_1(z_0\rightarrow z)\ U\  \gamma_2 (O)\  U\  \gamma_3(z \rightarrow z_0)
\end{equation}

where,

$\gamma_1(z \rightarrow z_0)$ : line segment from z to $z_0$;

$\gamma_2(O)$ : small circle centered at $z_0$;

$\gamma_3(z_0\rightarrow z)$ : line segment from $z_0$ to z.

\begin{figure}[h]
 \begin{center}
\includegraphics[width=0.6\linewidth]{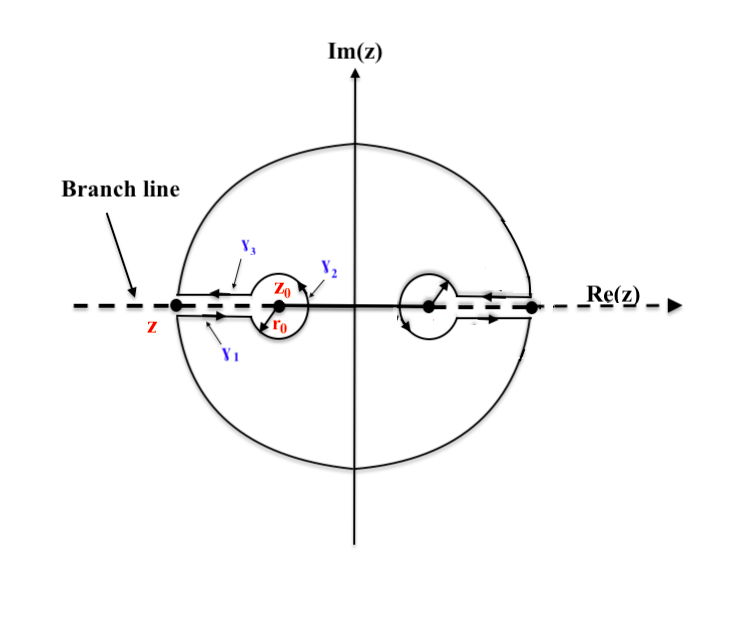}

\caption {contour of integration }
\end{center}

      \end{figure}

Then the Cauchy's integral formula becomes:

\begin{equation}
D^n f(z)=\frac{\Gamma(n +1)}{2\pi i} \int_{\gamma(z_0,z)} \frac{f(z)}{(z-z_0)^{n+1}}dz=I_{\gamma_1}+I_{\gamma_2}+I_{\gamma_3}
\end{equation}
$I_{\gamma_1}, I_{\gamma_2}, I_{\gamma_3}$ denote the integrals over the mentioned contours $\gamma_1, \gamma_2, \gamma_3$. Then, the line in which the branch occurs can be written as

\begin{equation}
\frac{1}{(z-z_0)^{\alpha+1}}= e^{(-\alpha-1)(ln |z-z_0|+i (\theta-\pi))}\ \ \ on\ \ \gamma_1
\end{equation}

\begin{equation}
\frac{1}{(z-z_0)^{\alpha+1}}=0\ \ \ on\ \ \gamma_2
\end{equation}

\begin{equation}
\frac{1}{(z-z_0)^{\alpha+1}}=e^{(-\alpha-1)(ln |z-z_0|+i (\theta+\pi))}\ \ \ on\ \ \gamma_3
\end{equation}

It should be noted that the integral tends to zero on $\gamma_2$ as the contour's radius $r_0$ goes to zero. Substituting the above equations inside of Eq. ($\ref{23}$) we obtain

\begin{equation}
^{\gamma(z,z^+)} D_{z-z_0}^\alpha f(z))=\small{\frac{ (e^{i \pi \alpha}+e^{-i \pi \alpha})  \Gamma(\alpha +1)}{2\pi i} \int_{z_0}^z \frac{f(z)}{(z-z_0)^{\alpha+1}}dz}
\end{equation}

or

\begin{equation}\label{32}
^{\gamma(z,z^+)} D_{z-z_0}^\alpha f(z)=\frac{ sin (\pi \alpha) \ \Gamma(\alpha +1)}{\pi } \int_{z_0}^z \frac{f(z)}{(z-z_0)^{\alpha+1}}dz.
\end{equation}

Notice that above equation is valid for all values of $\alpha$. By Weierstrass M-test we can show that an infinite series of functions converges. First we show that $\frac{f(z)}{z-z_0}$ is an infinite series:

\begin{equation}
^{\gamma(z,z^+)} D_{z-z_0}^\alpha f(z)=\frac{ sin (\pi \alpha) \ \Gamma(\alpha +1)}{\pi } \int_{z_0}^z \frac{f(z)}{(z-z_0)^{\alpha+1}}dz.
\end{equation}

\begin{equation}
=\frac{ sin (\pi \alpha) \ \Gamma(\alpha +1)}{\pi } \int_{z_0}^z \frac{f(z)}{(z-r_0)^\alpha}\ . \ \frac{1}{[1-(z_0-r_0)/(z-r_0)]^\alpha}
\end{equation}

\begin{equation}
=\frac{ sin (\pi \alpha) \ \Gamma(\alpha +1)}{\pi } \int_{z_0}^z \frac{f(z)}{(z-r_0)^\alpha}\ \Sigma_{n=0}^\infty (\frac{z_0-r_0}{z-r_0})^{n+\alpha}
\end{equation}

\begin{equation}
=\Sigma_{n=0}^\infty \frac{ sin (\pi \alpha) \ \Gamma(\alpha +1)}{\pi } \int_{z_0}^z \frac{(z_0-r_0)^{n+\alpha}}{(z-r_0)^{n+2\alpha}} f(z).
\end{equation}

As $\frac{f(z)}{(z-r_0)^\alpha}$ is bounded on $\gamma$, when the contour's radius goes to zero, by some positive number $M$, and $|\frac{z_0-r_0}{z-r_0}| \ \leq r<1$, then we have

\begin{equation}
|\frac{(z_0-r_0)^{n+\alpha}}{(z-r_0)^{n+2\alpha}} f(z)|\  \leq Mr^n
\end{equation}
this means that the series converges on $\gamma$.

Now we have an equation to be used for calculating our integral. Returning to our fractional integral

\begin{equation}
 \int_{0}^{2\pi}   \frac{ 1}{[1-\xi \ \cos(\psi)]^{3/2}}\ d\psi
\end{equation}
by writing the variable $\cos(\psi)$ in the complex plane as $\cos\psi=\frac{z+z^{-1}}{2}$, and replacing into the above equation, after some manipulations we obtain

\begin{equation}
 \int_{0}^{2\pi}   \frac{ 1}{[1-\xi \ \cos(\psi)]^{3/2}}\ d\psi=\small{\frac{2^{3/2}\pi}{\Gamma (3/2)} lim_{z\rightarrow z_0} D_{z}^{1/2}{(z-z_0)^{-3/2}  f(z)}}
\end{equation}

where $z_{01}=\frac{1+\sqrt{1-\xi^2}}{\xi}$ and $z_{02}=\frac{1-\sqrt{1-\xi^2}}{\xi}$ are the branch points of the integral in which the residues should be computed. We have then

\begin{equation}
 \int_{0}^{2\pi}   \frac{ 1}{[1-\xi \ \cos(\psi)]^{3/2}}\ d\psi=2\pi i\  \Sigma_{k=1}^n Res_f(z_k)
 \end{equation}
 
 \begin{equation}
=2\pi \   \Sigma Res \frac{f(z)}{z}
\end{equation}

\begin{equation}
=2\pi \   \Sigma Res \frac{1}{(1-\xi \frac{z+z^{-1}}{2})^{3/2}}\frac{1}{z}
\end{equation}

\begin{equation}
=2^{5/2}\pi \   \Sigma Res \frac{z^{1/2}}{\xi^{3/2}(-z^2+2(z/\xi)-1 )^{3/2}}.
\end{equation}

As there is a symmetry in the integral, it is not necessary to branch cut both the branch lines. For this reason we will take the interval [0,$\pi$] where one of the branch line is located,

\begin{equation}\label{81}
 \int_{0}^{\pi}   \frac{ 1}{[1-\xi \ \cos(\theta)]^{3/2}}\ d\psi=2\pi i\  \Sigma_{k=1}^n Res_f(z_k)
 \end{equation}
 
 \begin{equation}
=2^{5/2}\frac{\pi}{\Gamma (3/2)} [lim_{z\rightarrow \frac{1-\sqrt{1-\xi^2}}{\xi}} D_{z}^{1/2}(z-z_1)^{3/2} f(z)]
\end{equation}

\begin{equation}
=2^{5/2}\frac{\pi}{\Gamma (3/2)} [D_{z-z_1}^{1/2}\frac{z^{1/2}}{(z-(\frac{1-\sqrt{1-\xi^2}}{\xi}))^{3/2}}]
\end{equation}
where $D_{z}^{1/2}$ is the fractional derivative of the order $1/2$. Applying $D_{z}^{1/2}$ to the function we obtain,

\begin{equation}\label{82}
D_{z}^{1/2}\ (\frac{z^{1/2}}{(z-z_0)^{3/2}})
=\frac{\Gamma(3/2)}{(z-z_0)^{3/2}}+\frac{z\ \Gamma(5/2)}{2(z-z_0)^{5/2}}-\frac{z^2\ \Gamma(7/2)}{16(z-z_0)^{7/2}}+...
\end{equation}

Finally substituting the Eq. ($\ref{82}$) inside the $($\ref{81}$)$ we obtain:

\begin{equation}
 \int_{0}^{\pi}   \frac{ 1}{[1-\xi \ \cos(\theta)]^{3/2}}\ d\theta= 2^{5/2}\frac{\pi}{\Gamma (3/2)}
  [\frac{\Gamma(3/2)}{(z-z_0)^{3/2}}+\frac{z\ \Gamma(5/2)}{2(z-z_0)^{5/2}}-\frac{z^2\ \Gamma(7/2)}{16(z-z_0)^{7/2}}+...]
  \end{equation}
  
  \begin{equation}
=\frac{\pi}{\Gamma (3/2)} [\frac{\Gamma(3/2)}{(1-\xi^2)^{3/4}}+\frac{(1-\sqrt{1-\xi^2})\ \Gamma(5/2)}{4(1-\xi^2)^{5/4}}
-\frac{(1-\sqrt{1-\xi^2})^2\ \Gamma(7/2)}{64(1-\xi^2)^{7/4}}+...]
  \end{equation}

Where the integral's solution is a hypergeometric function and it can be written in a compact form as

\begin{equation}
I_1(\xi)=\int_{0}^{\pi}   \frac{ d\psi}{[1-\xi \ cos(\psi)]^{3/2}}=\frac{\pi}{(1+\xi)^{3/2}}\  _2F_1 (\frac{1}{2},\frac{3}{2};1;\frac{2\xi}{1+\xi}).
\end{equation}

By performing the same process using the so called Cauchy-like integral formula, Eq. ($\ref{19}$), we solved the integral ($\ref{19}$)

 \begin{equation}
I_2(\xi)=\int_{0}^{\pi}   \frac{ \cos(\psi)}{[1-\xi \ \cos(\psi)]^{3/2}}\ d\psi=\frac{\pi}{(1+\xi)^{3/2}}\  [\ _2F_1 (\frac{3}{2},\frac{3}{2};2;\frac{2\xi}{1+\xi})\ - \ _2F_1 (\frac{1}{2},\frac{3}{2};1;\frac{2\xi}{1+\xi})].
\end{equation}

\section{Conclusions}

In this paper, we used the Branch Cut method to change the non-local property of the fractional derivative to a local property in order to be used the Cauchy's integral formula. The new method can be called the Cauchy-like integral formula for fractional integrals. By Weierstrass M-test we have shown that the integral converges. This method helps to solve important problems in which there is branch line. We applied the method to solved the axial and radial magnetic field components for a coil of negligible thickness with a stationary electric current. At the end, we have shown that the solutions can be expressed by means of the hypergeometric functions $_2F_1(a,b;c;z)$.





\end{document}